\documentclass[12pt]{article}

\def\Xd{X^{(d)}}
\def\Xi{X^{(\infty)}}
\def\Zd{Z^{(d)}}
\def\TV#1{\|#1\|_{TV}}
\def\KR#1{\|#1\|_{KR}}
\def\lip{{\rm Lip}_1^1}
\def\dist{{\rm dist}}

\def\X{{\cal X}}

\def\bx{{\bf x}}

\def\bY{{\bf Y}}
\def\bZ{{\bf Z}}

\usepackage{jsr}
\usepackage{amssymb}

\def\proof{\medskip\noindent\bf Proof. \rm \ }

\begin{document}

\baselineskip 24pt

\begin{center}
{\titlefont Complexity Bounds for MCMC via Diffusion Limits}
\\
by \\
\gareth\ and \jeff\ \\
(August 15, 2014.) \\
\end{center}

\abstract{
We connect known results about diffusion limits of Markov chain
Monte Carlo (MCMC) algorithms to
the Computer Science notion of algorithm complexity.  Our main result
states that any diffusion limit of a Markov process implies a
corresponding complexity bound (in an appropriate metric).  We then combine
this result with previously-known MCMC diffusion limit results
to prove that under
appropriate assumptions, the Random-Walk Metropolis (RWM) algorithm in $d$
dimensions takes $O(d)$ iterations to converge to stationarity, while
the Metropolis-Adjusted Langevin Algorithm (MALA)
takes $O(d^{1/3})$ iterations to converge to stationarity.
}

\sect{Introduction}

In the computer science literature, algorithms are often analysed
in terms of ``complexity'' bounds.  In the Markov chain Monte Carlo
(MCMC) literature, algorithms are sometimes understood in terms of
diffusion limits.  The purpose of this note is to connect these two
approaches, and in particular to show that diffusion limits sometimes
imply complexity bounds.

Complexity results in computer science go back at least to Cobham
(1964), and took on greater focus with the pioneering {\it NP-complete}
work of Cook (1971).  In the Markov chain context, computer scientists
have been bounding convergence times of Markov chain algorithms since
at least Jerrum and Sinclair (1989), focusing largely on spectral
gap bounds for Markov chains on finite state spaces.  More recently,
attention has turned to bounding spectral gaps of modern Markov chain
algorithms on general (e.g.\ uncountable) state spaces, again primarily
via spectral gaps (e.g.\ Woodard et al., 2009a, 2009b).  These bounds
often focus on the order of the convergence time in terms of some
parameter such as the dimension $d$ of the corresponding state space.

Meanwhile, in statistics, MCMC algorithms are extremely widely used
and studied (see e.g.\ Brooks et al., 2011, and the many references
therein), and their running times are an extremely
important practical issue.  They have been studied from a variety
of perspectives, including directly bounding the convergence in
total variation distance (see e.g.\ Rosenthal, 1995b, 1996, 2002;
Jones and Hobert, 2001, 2004; and references therein), convergence
``diagnostics'' via statistical analysis of the Markov chain output
(e.g.\ Gelman and Rubin, 1992), and most notably by proving weak
convergence limits of sped-up versions of the algorithms to diffusion
limits (e.g.\ Roberts et al., 1997; Roberts and Rosenthal, 1998).

The MCMC direct total variation bounds are sometimes presented in terms
of the convergence order (e.g.\ see Rosenthal, 1995a, for order bounds
for a Gibbs sampler for a variance components model).  In addition,
the MCMC diffusion limits often involve speeding up the original
algorithm by a certain order, and then proving weak convergence to a
fixed process which converges in $O(1)$ iterations, thus giving them
the flavour of complexity order bounds too.  However, the MCMC
results are typically not stated precisely in terms of convergence
time complexity results, and (perhaps because of this) they are often
overlooked by the computer science complexity community.

In this paper, we attempt to connect these two streams of
Markov chain convergence time bounds.  In particular, we establish
(Theorem~\ref{mainthm}) that results about diffusion limits do directly
imply corresponding complexity bounds (using an appropriate
convergence metric as described below).
We then apply our theorem to previous results about diffusion limits
of MCMC algorithms (Section~\ref{sec-MCMC}),
to establish running time complexity order bounds
for such MCMC algorithms as the Random-Walk Metropolis algorithm (Theorem~\ref{RWMthm}) and
the Metropolis-adjusted Langevin algorithm (Theorem~\ref{MALAthm}).

\sect{Assumptions and Main Result}

Let $(\X,\F,\rho)$ be a general
measurable metric space, i.e.\ a non-empty (and possibly uncountable)
set $\X$ endowed with a
metric $\rho$ which induces a Borel $\sigma$-algebra $\F$
of measurable subsets.
We wish to bound the convergence of a stochastic
process $\{X_t\}$ on $(\X,\F)$ to its stationary
probability distribution $\pi$.  To measure the distance to
stationarity, on finite state spaces one often (see e.g.\ Aldous and Fill,
2002, Section 2.4.1) uses the total variation distance defined by
$$
\|\L_x(X_t) - \pi\|_{TV}
\ := \
\sup_{|f| \le 1}
\Big| \E_x[f(X_t)] - \pi(f) \Big|
$$
where the supremum is taken over all measurable functions
$f:\X\to\IR$ with $|f(x)| \le 1$ for all $x\in\X$.
Here $\L_x(X_t)$ is the law of $X_t$ conditional on starting at
$X_0=x$, and $\E_x[f(X_t)]$ is the expected value of $f$ with respect
to this law, and $\pi(f) = \int f(x) \, \pi(dx)$ is the expected
value of $f$ with respect to $\pi$.

This total variation distance can also be used on general state spaces
in many instances (see e.g.\ Rosenthal, 1995).  However, it is not appropriate
for bounding the weak convergence which arises in the diffusion context,
since it may not go to zero for processes which converge only weakly to
stationarity,
so we do not use it here.  Instead, we let
$$
\lip \ = \ \Big\{ f : \X\to\IR, \ |f(x)-f(y)| \le \rho(x,y) \ \forall
x,y\in\X, \ |f| \le 1 \Big\}
$$
be the set of all
functions from $\X$ to $\IR$ with Lipschitz constant $\le 1$ {\it and}
with $|f(x)| \le 1$ for all $x\in\X$, and use the distance function
$$
\|\L_x(X_t) - \pi\|_{KR}
\ := \
\sup_{f \in \lip}
\Big| \E_x[f(X_t)] - \pi(f) \Big|
\, .
$$
(Here ``KR'' stands for
``Kantorovich-Rubinstein''; see the proof of
Proposition~\ref{weakconvprop} below.)
The distance $\KR{\cdots}$ is similar to, but more restrictive than,
the total variation distance, and we will see below
(Proposition~\ref{weakconvprop}) that it metrises
weak convergence and so is appropriate for our purposes.

We also note that many approaches to stationary instead directly
bound the spectral gap of the corresponding Markov operator (e.g.\
Woodard et al., 2009b).  However, on general state spaces,
the spectral gap is zero for Markov chains which
are not ``geometrically ergodic'' (see e.g.\ Theorem~2 of Roberts and
Rosenthal, 1997).  Furthermore, many MCMC algorithms are not
geometrically ergodic (e.g.\ the Random-Walk Metropolis algorithm
on target distributions with heavier-than-exponential tails, see
Theorem~3.3 of Mengersen and Tweedie, 1996).
They also are often not reversible, which makes spectral gaps
harder to study or interpret.
For these reasons, we do not wish to restrict attention to
spectral gaps, which is another reason that we use the metric $\KR{\cdots}$.

A related issue is what initial states $X_0$ should be considered.  On
finite state spaces, one often (e.g.\ Jerrum and Sinclair, 1989,
Section~2) considers the worst case, by taking
supremum over all initial states $x$, i.e.\ uses something like
$\sup_{x\in\X} \|\L_x(X_t) - \pi\|_{TV}$.  But this supremum is also
frequently inappropriate on general state spaces.  For instance,
if $\X$ is unbounded,
then as $t$ increases one can start from worse and worse states $X_0$ so
that the supremum will
never go to~0.  Instead, we need to specify more precisely which initial
state(s) $X_0$ to consider.  As a concrete choice, we will
take the $\pi$-average of the distances to stationarity
from all initial states $X_0$ in $\X$.  That is, for any Markov chain $\{X_t\}$
on $(\X,\F)$ with stationary distribution $\pi$, we measure the distance
to stationarity at time $t$ by the distance function
$$
\E_{X_0\sim\pi} \|\L_{X_0}(X_t) - \pi\|_{KR}
\ := \
\int_{x\in\X} \pi(dx) \ \|\L_x(X_t) - \pi\|_{KR}
\, .
$$

Using this distance function, we can state our main result:

\pnew{Theorem}
\plabel{mainthm}
Let $\Xd = \{\Xd_t\}_{t \ge 0}$ be a stochastic
process on $(\X,\F,\rho)$,
for each $d\in\IN$, which
converges weakly in the
Skorokhod topology as $d\to\infty$
to another stochastic process $\Xi = \{\Xi_t\}_{t \ge 0}$,
i.e.\ $\Xd_t \wto \Xi_t$ for each fixed $t \ge 0$.
Assume these processes all have the same stationary
probability distribution $\pi$, and that $\Xi$ converges
(either weakly or in total variation distance) to $\pi$.
Then for any $\epsilon>0$,
there are $D<\infty$ and $T<\infty$ such that
$$
\E_{\Xd_0\sim\pi} \KR{\L_{\Xd_0}(\Xd_t) - \pi} \ < \ \epsilon \, ,
\quad t \ge T,
\quad d \ge D \, .
$$

Theorem~\ref{mainthm} may be summarised as saying that if a sequence
$\{\Xd\}$ of Markov processes converges weakly to a limiting ergodic
process, then we can bound the convergence of the sequence of processes
uniformly over all sufficiently large $d$, i.e.\ the processes converge
in $O(1)$ iterations with respect to $d$.  We will next apply this
result to previously known diffusion limits of common MCMC algorithms.

\sect{Application to MCMC}
\labels{sec-MCMC}

Our primarily interest is in the use of Theorem~\ref{mainthm} to bound
the complexity of MCMC algorithms.  We begin with the most popular
MCMC algorithm, the Random-Walk Metropolis (RWM) algorithm.  This algorithm
proceeds, given a positive target probability density $\pi_d$ on the state
space $\IR^d$, by running a Markov chain $\{\bZ^d_n\}_{n=0}^\infty$ as
follows.  Given the value $\bZ^d_n$, a proposed new state
$\bY^d_{n+1} \sim MVN(\bZ^d_n, \, \sigma^2_d)$ is chosen from a
multivariate normal distribution centered at $\bZ^d_n$, and then with
probability $\min[1, \ \pi(\bY^d_{n+1})/\pi(\bZ^d_n)]$ the proposal is
accepted and $\bZ^d_{n+1} = \bY^d_{n+1}$, otherwise with the remaining
probability the proposal is rejected and $\bZ^d_{n+1} = \bZ^d_{n}$.
This algorithm is easily seen to be irreducible and aperiodic and
to leave $\pi$ stationary, so it will converge asymptotically to
$\pi$.  The question then becomes how quickly it will converge, and
what choice of proposal variance $\sigma^2_d$ is optimal.

In this context,
Roberts et al.\ (1997) proved the remarkable result that $U^d \wto U$
as $d\to\infty$, where
$U^d_t = \bZ^d_{\lfloor dt \rfloor, \, 1}$ is the first coordinate of the
RWM algorithm
sped up by a factor of $d$, and $U$ is a limiting ergodic Langevin
diffusion, and $\wto$ indicates weak convergence in the usual
Skorokhod topology.
They proved this result under certain strong technical assumptions, namely
that $\pi_d$ takes on the special product form
$\pi_d(\bx) = \prod_{i=1}^d h(x_i)$ for some fixed function
$h:\IR\to(0,\infty)$ with $h'/h$ Lipschitz continuous,
and $\int [h'(x)/h(x)]^8 h(x) dx < \infty$,
and $\int [h''(x)/h(x)]^4 h(x) dx < \infty$.
They also assumed the other coordinates 2 through $d$ of
the process $\bZ^d$ are in stationarity, and
that $\sigma^2_d = \ell^2/(d-1)$ for some fixed $\ell>0$.

This theorem of Roberts et al.\ (1997) allowed them to study the
limiting diffusion $U$ as a function of the proposal variance
parameter $\ell$, and optimise it to prove that the algorithm
converges fastest when its asymptotic acceptance rate is equal to
0.234\ldots\ (see also Roberts and Rosenthal, 2001).  Furthermore, since
their process $U^d$ involved speeding up the original algorithm by a
factor of $d$, their results seemed to imply that RWM required $O(d)$
iterations to converge.  However, a precise statement of such a
complexity bound was not provided.

In light of Theorem~\ref{mainthm} above, we are now able to use the
diffusion limit of Roberts et al.\ (1997) to give an actual complexity
bound on the RWM algorithm.  Indeed, applying Theorem~\ref{mainthm}
to their limit immediately yields:

\pnew{Theorem}  
\plabel{RWMthm}
Let $\Zd$ be a RWM algorithm on a product density in $d$ dimensions
satisfying the technical assumptions of Roberts et al.\ (1997).
Then for any $\epsilon>0$,
there is $D<\infty$ and $T<\infty$ such that
$$
\E_{\Zd_0\sim\pi} \KR{\L_{\Zd_0}(\Zd_{\lfloor dt \rfloor,1}) - h}
\ < \ \epsilon \, ,
\quad t \ge T,
\quad d \ge D \, .
$$
Hence, the RWM algorithm takes $O(d)$ iterations to converge
to within $\epsilon$ of stationarity in any one coordinate.

We believe this to be the first precise general result about the
convergence order of the RWM algorithm.  Of course, it requires
the strong technical assumptions of Roberts et al.\ (1997), but it
still applies to a fairly general collection of densities on $\IR^d$.
Furthermore, it appears empirically (see e.g.\ Roberts and Rosenthal,
2001) that even when RWM algorithms do not satisfy the technical
assumptions they still exhibit similar limiting behaviour.

Another MCMC diffusion limit concerns the Metropolis-Adjusted Langevin
Algorithm (MALA).  This algorithm is similar to the above Random-Walk
Metropolis algorithm, except that now the proposal state
$\bY^d_{n+1} \sim MVN(\bZ^d_n + \half \sigma_d^2 \grad \log
\pi_d(Z^d_n), \, \sigma^2_d)$ is chosen from a
multivariate normal distribution centered at $\bZ^d_n + \half
\sigma_d^2 \grad \log
\pi_d(Z^d_n), \, \sigma^2_d)$ (to better approximate $\pi$), and the
above acceptance probability is modified by the ratio of the
corresponding proposal normal distributions.
In this context, Roberts and Rosenthal (1999) proved
that $U^d \wto U$, where $U^d_t = \bZ^d_{\lfloor d^{1/3} t \rfloor, \, 1}$
is the first coordinate of the MALA algorithm
sped up by a factor of $d^{1/3}$, and $U$ is again a limiting
ergodic Langevin diffusion.
This result again required strong technical assumptions, this time that 
$\pi_d(\bx) = \prod_{i=1}^d h(x_i)$ for some fixed function
$h:\IR\to(0,\infty)$ with polynomially-bounded log-derivatives of all
orders, and finite moments of all orders, with $h'/h$ Lipschitz continuous.
They also assumed that coordinates 2 through $d$ of
$\bZ^d$ are again in stationarity, and
that $\sigma^2_d = \ell^2 \, d^{-1/3}$ for some fixed $\ell>0$.

This theorem of Roberts and Rosenthal (1999) allowed them to optimise
the limiting diffusion $U$ as a function of $\ell$, and to prove that
the algorithm converges fastest when its asymptotic acceptance rate
is equal to 0.574\ldots.  Also, since their process $U^d$ involved
speeding up the original algorithm by a factor of $d^{1/3}$, their
results seemed to imply that MALA required $O(d^{1/3})$ iterations
to converge.  Once again, we can use Theorem~\ref{mainthm} above to
obtain a more formal complexity bound:

\pnew{Theorem}  
\plabel{MALAthm}
Let $\Zd$ be a MALA algorithm on a product density in $d$ dimensions
satisfying the technical assumptions of Roberts and Rosenthal (1999).
Then for any $\epsilon>0$,
there is $D<\infty$ and $T<\infty$ such that
$$
\E_{\Zd_0\sim\pi} \KR{\L_{\Zd_0}(\Zd_{\lfloor d^{1/3} t \rfloor,1}) - h}
\ < \ \epsilon \, ,
\quad t \ge T,
\quad d \ge D \, .
$$
Hence, the MALA algorithm takes $O(d^{1/3})$ iterations to converge
to within $\epsilon$ of stationarity in any one coordinate.

Finally, we note that a number of other diffusion limits have been
proven for MCMC algorithms in other contexts.  For example, B\'edard
(2007, 2008) and Sherlock and Roberts (2009) have extended the
original RWM diffusion limit to more general target distributions;
Roberts (1998) and Neal and Roberts (2006, 2008, 2011)
and Jourdain et al.\ (2013a, 2013b)
have extended it to other related cases; and
Neal et al.\ (2012) have established diffusion limits for RWM
algorithms on discontinuous target densities.  Each of these diffusion
limit results could also be combined with Theorem~\ref{mainthm} above
to yield complexity order bounds in new contexts.

\sect{Proof of Theorem~\ref{mainthm}}

In this section, we prove Theorem~\ref{mainthm}.  Along the
way, we establish that $\KR{\cdots}$ metrises weak convergence
(Proposition~\ref{weakconvprop}), and that $\E_{X_0\sim\pi}
\KR{\L_{X_0}(\Xd_t) - \pi}$ is a non-increasing function of $t$
(Lemma~\ref{monexplemma}).  We first establish that $\KR{\cdots}$
is a norm:


\pnew{Lemma}
\plabel{normlemma}
Let $S$ be any non-empty collection of functionals $\X\to\IR$
which is symmetric (i.e.\ if $f\in S$ then $-f\in S$).
Let $\|\mu\| = \sup\limits_{f\in S} \mu(f)$.
Then $\|\ldots\|$ is a (possibly infinite) norm function
on the set of all signed measures on $(\X,\F)$.
In particular, $\KR{\cdots}$ is a norm.

\proof It is immediate that $\|0\|=0$,
and that $\|a \, \mu\| = a \, \|\mu\|$ for $a > 0$.
The symmetry of $S$ implies that
$\|-\mu\| = \|\mu\|$.
Finally, for the triangle inequality, we check that
$$
\|\mu+\nu\| \ = \ \sup_{f \in S} \Big( \mu(f) + \nu(f) \Big)
\ \le \ \Big( \sup_{f \in S} \mu(f) \Big) + \Big( \sup_{f \in S} \nu(f) \Big)
\ = \ \|\mu\| + \|\nu\|
\, .
$$
Hence, $\|\ldots\|$ is a norm.
The claim about $\KR{\cdots}$ then follows by taking $S = \lip$.
\qed

We next show that
truncating the metric $\rho$ does not change $\lip$:

\pnew{Lemma}
\plabel{trunclemma}
Let $\rho^* = \min(2, \, \rho)$.  Then
$$
\lip \ = \ \{f : \X\to\IR, \ |f(x)-f(y)| \le \rho^*(x,y) \ \forall
x,y\in\X, \ |f| \le 1\}
\, .
$$

\proof This is immediate since we always have $|f(x)-f(y)| \le 2$
for $f\in\lip$.
\qed

\pnew{Proposition}
\plabel{weakconvprop}
The metric
$\Delta(\mu,\nu) := \KR{\mu-\nu}$ metrises weak convergence of
probability measures on $(\X,\F,\rho)$.  That is,
if $\{\mu_t\}$ and $\mu$ are probability measures on $(\X,\F,\rho)$,
then $\{\mu_t\} \wto \mu$ if and only if $\lim_{t\to\infty}
\Delta(\mu_t,\mu) = 0$.

\proof
Let $\rho^*$ be as in Lemma~\ref{trunclemma}.
We first note that since $\rho$ and $\rho^*$ agree for distances $\le 2$,
they give rise to precisely the same open subsets.  Therefore,
$(\X,\rho^*)$ induces the same Borel $\sigma$-algebra $\F$ that $(\X,\rho)$
does, and thus gives rise to the same Skorokhod topology.  Hence, weak
convergence on $(\X,\F,\rho)$ is precisely equivalent to
weak convergence on $(\X,\F,\rho^*)$.  Furthermore, by
Lemma~\ref{trunclemma}, the metric $\KR{\cdots}$ is the same on
$(\X,\F,\rho^*)$ as on $(\X,\F,\rho)$.
Hence, it suffices to prove the result 
on the truncated space $(\X,\F,\rho^*)$.

Now, since $(\X,\F,\rho^*)$ is a bounded metric space,
it is known (see e.g.\ Givens and Shortt, 1984, Proposition~4)
that weak convergence on $(\X,\F,\rho^*)$ is metrised by the
Wasserstein metric $W_1$ on $(\X,\rho^*)$, defined by
$$
W_1(\mu,\nu) \ := \ \inf \E[\rho(X,Y)]
$$
where the infimum is taken over all pairs $(X,Y)$ of random variables
on $(\X,\F)$ such that $\L(X)=\mu$ and $\L(Y)=\nu$.
On the other hand, again since $(\X,\F,\rho^*)$ is a bounded metric space,
it is known
(Kantorovich and Rubinstein, 1958; see e.g.\ Givens and Shortt, 1984, p.~233)
that for probability measures $\mu$ and $\nu$ on $(\X,\F,\rho^*)$,
the Wasserstein metric $W_1(\mu,\nu)$ is precisely
equal to $\KR{\mu-\nu}$.
Combining these two facts, the result follows
for $(\X,\F,\rho^*)$, and hence also for $(\X,\F,\rho)$.
\qed

\pnew{Lemma}
\plabel{Xilemma}
If $\Xi$ converges to $\pi$, either weakly or in total variation
distance, then for all $\epsilon>0$ there is $T<\infty$ such that
$\KR{\L_x(\Xi_T) - \pi} \le \epsilon/2$ for all $t \ge T$.

\proof
If the convergence is weak, then this follows from
Proposition~\ref{weakconvprop}.
If the convergence is
in total variation distance, then
this still follows since $\KR{\ldots} \le \TV{\ldots}$.
\qed

\pnew{Proposition}
\plabel{Tconvprop}
Under the assumptions of Theorem~\ref{mainthm}, for any $x\in\X$ and
$\epsilon>0$,
there is $D<\infty$ and $T<\infty$ such that
$$
\KR{\L_x(\Xd_T) - \pi} \ < \ \epsilon \, ,
\quad d \ge D \, .
$$

\proof
Using Lemma~\ref{normlemma}, we have by the triangle inequality that
$$
\KR{\L_x(\Xd_t) - \pi}
\ \, \le \ \, \KR{\L_x(\Xd_t) - \L_x(\Xi_t)}
\ + \ \KR{\L_x(\Xi_t) - \pi}
\, .
$$
By Lemma~\ref{Xilemma},
there is $T<\infty$ such that
$\KR{\L_x(\Xi_T) - \pi} \le \epsilon/2$.
Then, since $\Xd_T$ converges weakly to $\Xi_T$,
by Proposition~\ref{weakconvprop}
there is $D<\infty$ such that for all $d \ge D$,
$\KR{\L_x(\Xd_T) - \L_x(\Xi_T)} < \epsilon/2$.
The result follows.
\qed

\rnew{Remark}
\plabel{Sremark}
If the weak
convergence of $\Xd$ to $\Xi$ is assumed to be uniform over 
bounded time intervals, then we can strengthen
Proposition~\ref{Tconvprop} to say that for any $x\in\X$ and
$\epsilon>0$ and $S<\infty$, there are
$D<\infty$ and $T<\infty$ such that
$\KR{\L_x(\Xd_t) - \pi} < \epsilon$ for all
$t \in [T, \, T+S]$.

\pnew{Corollary}
\plabel{Tconvcor}
Under the assumptions of Theorem~\ref{mainthm}, for any $\epsilon>0$,
there is $D<\infty$ and $T<\infty$ such that
$$
\E_{X_0\sim\pi} \KR{\L_{X_0}(\Xd_T) - \pi} \ < \ \epsilon \, ,
\quad d \ge D \, .
$$

\proof
We first let
$$
A_m \ = \ \{x\in\X : \KR{\L_x(\Xd_t) - \pi} \, < \, \epsilon/2
\ \ \forall t \ge m\}
\, .
$$
Then $A_{m+1} \subseteq A_m$ by inspection, and
$\bigcup_m A_m = \X$ by Lemma~\ref{Xilemma}.
Hence, by continuity of probabilities (see e.g.\ Proposition~3.3.1
of Rosenthal, 2000), $\lim_{m\to\infty} \pi(A_m) = 1$.
We can therefore find $T<\infty$ such that $\pi(A_T) \ge 1 - (\epsilon/8)$.

Next, for this fixed $T$, let
$$
B_m \ = \ \{x\in\X : \KR{\L_x(\Xd_T) - \pi} \, < \, \epsilon/2
\ \ \forall d \ge m\}
\, .
$$
Then $B_{m+1} \subseteq B_m$ by inspection, and
$\bigcup_m B_m = \X$ by Proposition~\ref{Tconvprop},
so again by continuity of probabilities
we can find $D\in\IN$ such that $\pi(B_D) \ge 1 - (\epsilon/8)$.

We then compute that for this fixed $T$ and $D$, and for any $d \ge D$,
$$
\E_{X_0\sim\pi} \KR{\L_{X_0}(\Xd_T) - \pi}
\qquad \qquad \qquad \qquad
\qquad \qquad \qquad \qquad
\qquad \qquad \qquad \qquad
$$
$$
\ = \ \E_{X_0\sim\pi} \Big( \one_{{X_0} \in A_T \cap B_D}
	\, \KR{\L_{X_0}(\Xd_T) - \pi} \Big)
\, + \, \E_{X_0\sim\pi} \Big( \one_{{X_0} \not\in A_T \cap B_D}
	\, \KR{\L_{X_0}(\Xd_T) - \pi} \Big)
$$
$$
\ \le \ (\epsilon/2) + [(\epsilon/8)+(\epsilon/8)] \times 2
\ = \ \epsilon \, ,
$$
where we have used the fact that by definition we always have
$\KR{\L_x(\Xd_T) - \pi} \le 2$ for any $x$ and $d$.
This gives the result.
\qed



Corollary~\ref{Tconvcor} is nearly what we need to prove
Theorem~\ref{mainthm}.  However, for Theorem~\ref{mainthm} we want the
convergence to
be within $\epsilon$ for {\it all} $t \ge T$, not just for one fixed $T$
(nor just for all $t$ in some bounded time interval,
cf.\ Remark~\ref{Sremark}).
Unfortunately, $\KR{\L_x(\Xd_t) - \pi}$
might not be a non-increasing function of $t$
(unlike $\TV{\L_x(\Xd_t) - \pi}$, which always is, see e.g.\
Proposition~3(c) of Roberts and Rosenthal, 2004).
On the other hand, fortunately the quantity
$\E_{X_0\sim\pi} \KR{\L_{X_0}(\Xd_t) - \pi}$ is indeed non-increasing:

\pnew{Lemma}
\plabel{monexplemma}
Let $\|\ldots\|$ be any norm function on signed measures on $(\X,\F)$.
Let $P^t(x,\cdot)$ be the transition probabilities for a Markov chain
on $(\X,\F)$ with stationary probability distribution $\pi$.
Let $\dist(t) = \E_{X_0\sim\pi}\|P^t(X_0,\cdot)-\pi\|$.
Then $\dist(t)$ is a non-increasing function of $t$.
In particular, in the context of Theorem~\ref{mainthm},
$\E_{X_0\sim\pi} \KR{\L_{X_0}(\Xd_t) - \pi}$ is a non-increasing function
of $t$.

\proof
We compute by stationarity that for $s,t > 0$,
$$
\dist(s+t) \ = \ \E_{X_0\sim\pi}\|P^{s+t}(X_0,\cdot)-\pi\|
$$
$$
\ = \ \E_{X_0\sim\pi}\left\|\int_{y\in\X} P^s(X_0,dy)
	\ P^{t}(y,\cdot)-\pi\right\|
$$
$$
\ \le \ \E_{X_0\sim\pi}\int_{y\in\X} P^s(X_0,dy)
\, \Big\|P^{t}(y,\cdot)-\pi\Big\|
$$
$$
\ = \ \E_{Y_0\sim\pi} \|P^{t}(Y_0,\cdot)-\pi\|
\ = \ \dist(t)
\, ,
$$
thus proving the first claim.  The claim about
$\E_{x\sim\pi} \KR{\L_x(\Xd_t) - \pi}$ then follows by
Lemma~\ref{normlemma} upon setting
$P^t(x,A) = \P[\Xd_t \in A \, | \, \Xd_0 = x]$.
\qed


\medskip
Theorem~\ref{mainthm} then follows by combining Corollary~\ref{Tconvcor}
and Lemma~\ref{monexplemma}.

\medskip
\ack We thank Dawn Woodard and Alexandre Thiery for helpful
discussions of these matters.

\bigskip
\section*{References}
\frenchspacing

D.~Aldous and J.A.~Fill (2002),
Reversible Markov Chains and Random Walks on Graphs.
Unfinished monograph, available at:
http://www.stat.berkeley.edu/$\sim$aldous/RWG/book.html

M.~B\'edard (2007),
Weak Convergence of Metropolis Algorithms for
Non-iid Target Distributions.
Ann. Appl. Prob. {\bf 17}, 1222--1244.

M.~B\'edard (2008),
Optimal Acceptance Rates for Metropolis Algorithms:
Moving Beyond 0.234.
Stoch. Proc. Appl. {\bf 118}, 2198--2222.

S.~Brooks, A.~Gelman, G.L.~Jones, and X.-L.~Meng, eds.\ (2011),
Handbook of Markov chain Monte Carlo.
Chapman \& Hall / CRC Press.

A.~Cobham (1964),
The intrinsic computational difficulty of functions.
In Proceedings of the 1964
International Congress for Logic, Methodology, and Philosophy of
Science, Y.~Bar-Hille, ed.,
Elsevier/North-Holland, Amsterdam, 24-­30.

S.~Cook (1971), The complexity of theorem-proving procedures.
Third Annual ACM Symposium on Theory of Computing, 151-­158.

A. Gelman and D.B. Rubin (1992), Inference from iterative simulation
using multiple sequences.  Stat. Sci. {\bf 7(4)}, 457--472.

C.R.~Givens and R.M.~Shortt (1984),
A class of Wasserstein metrics for probability distributions.
Michigan Math.\ J.\ {\bf 31(2)}, 231--240.

G.L. Jones and J.P. Hobert (2001), Honest exploration of intractable
probability distributions via Markov chain Monte Carlo. Statistical
Science {\bf 16}, 312--334.

G.L. Jones and J.P. Hobert (2004), Sufficient burn-in for Gibbs samplers
for a hierarchical random effects model.  Ann. Stat. {\bf 32}, 784--817.

B.~Jourdain, T.~Leli\`evre, and B.~Miasojedow (2013a),
Optimal scaling for the transient phase
of Metropolis Hastings algorithms: the longtime behavior.
Bernoulli, to appear.

B.~Jourdain, T.~Leli\`evre, and B.~Miasojedow (2013b),
Optimal scaling for the transient phase
of Metropolis Hastings algorithms: the mean-field limit.
Preprint.

L.~Kantorovich and G.~Rubinstein (1958),
On a space of completely additive functions.  (In Russian.)
Vestnik Leningrad.\ Univ.\ {\bf 13}, 52--59.

K.L.~Mengersen and R.L.~Tweedie (1996),
Rates of convergence of the Hastings and Metropolis algorithms.
Ann. Stat. {\bf 24}, 101--121.

P.~Neal and G.O.~Roberts (2006),
Optimal Scaling for partially updating MCMC algorithms.
Ann. Appl. Prob. {\bf 16(2)}, 474--515.

P.~Neal and G.O.~Roberts (2008),
Optimal Scaling for Random Walk Metropolis on spherically
constrained target densities.
Meth. Comput. Appl. Prob. {\bf 10}, 277-­297.

P.~Neal and G.O.~Roberts (2011),
Optimal scaling of random walk Metropolis algorithms with non-Gaussian
proposals.
Meth. Comput. Appl. Prob. {\bf 13(3)}, 583--601.

P.~Neal, G.O.~Roberts, and W.K.~Yuen (2012),
Optimal scaling of random walk Metropolis algorithms with
discontinuous target densities.
Ann. Appl. Prob. {\bf 22(5)}, 1880--1927.

G.O.~Roberts (1998),
Optimal Metropolis algorithms on the hypercube.
Stochastics {\bf 62}, 275--284.

G.O.~Roberts, A.~Gelman, and W.R.~Gilks (1997),
Weak convergence and optimal scaling of random walk Metropolis algorithms.
Ann.\ Appl.\ Prob.\ {\bf 7}, 110--120.

G.O.~Roberts and J.S.~Rosenthal (1997),
Geometric ergodicity and hybrid Markov chains.
Elec. Comm. Prob. {\bf 2(2)}, 13--25.

G.O.~Roberts and J.S.~Rosenthal (1998),
Optimal scaling of discrete approximations to Langevin diffusions.
J. Roy. Stat. Soc. B {\bf 60}, 255--268.

G.O.~Roberts and J.S.~Rosenthal (2001),
Optimal scaling for various Metropolis-Hastings algorithms.
Stat. Sci. {\bf 16}, 351--367.

G.O.~Roberts and J.S.~Rosenthal (2004),
General state space Markov chains and MCMC algorithms.
Prob.\ Surv.\ {\bf 1}, 20--71.

J.S. Rosenthal (1995a), Rates of convergence for Gibbs sampler for
variance components models.  Ann. Stat. {\bf 23}, 740--761.

J.S. Rosenthal (1995b), Minorization conditions and convergence rates
for Markov chain Monte Carlo. J. Amer. Stat. Assoc. {\bf 90}, 558--566.

J.S. Rosenthal (1996), Convergence of Gibbs sampler for a model
related to James-Stein estimators.  Stat. and Comput. {\bf 6}, 269--275.

J.S. Rosenthal (2000),
A first look at rigorous probability theory.
World Scientific Publishing, Singapore.

J.S. Rosenthal (2002),
Quantitative convergence rates of Markov chains: A simple account.
Elec. Comm. Prob. {\bf 7}, No.~13, 123--128.

C.~Sherlock and G.O.~Roberts (2009),
Optimal scaling of the random walk Metropolis
on elliptically symmetric unimodal targets.
Bernoulli {\bf 15(3)}, 774--798.

D.B. Woodard, S.C. Schmidler, M.L. Huber (2009a),
Sufficient conditions for torpid mixing of parallel and simulated tempering.
Elec.\ J.\ Prob.\ {\bf 14}, 780--804.

D.B. Woodard, S.C. Schmidler, M.L. Huber (2009b),
Conditions for rapid mixing of parallel and simulated tempering
on multimodal distributions.
Ann.\ Appl.\ Prob.\ {\bf 19}, 617--640.

\end{document}